\documentclass[12pt]{article}

\usepackage{amssymb}
\usepackage{amsmath}
\usepackage{amsthm}
\usepackage[mathscr]{eucal}
\usepackage{cite}

\begin{document}

\renewcommand{\citeleft}{{\rm [}}
\renewcommand{\citeright}{{\rm ]}}
\renewcommand{\citepunct}{{\rm,\ }}
\renewcommand{\citemid}{{\rm,\ }}

\newcounter{abschnitt}
\newtheorem{koro}{Corollary}[abschnitt]
\newtheorem{theorem}[koro]{Theorem}
\newtheorem{satz}{Theorem}
\newtheorem{lem}[koro]{Lemma}
\newtheorem{prop}[koro]{Proposition}

\newcounter{saveeqn}
\newcommand{\alpheqn}{\setcounter{saveeqn}{\value{abschnitt}}
\renewcommand{\theequation}{\mbox{\arabic{saveeqn}.\arabic{equation}}}}
\newcommand{\reseteqn}{\setcounter{equation}{0}
\renewcommand{\theequation}{\arabic{equation}}}

\newcommand{\vol}{\mbox{vol}}

\newcommand{\R}{\mathbb{R}}

\renewcommand{\P}{ P_{e_i^{\perp}} }
\newcommand{\Pu}{ P_{u_i^{\perp}} }

\newcommand{\Eins}{\mathbbm{1}}


\sloppy

\phantom{a}

\vspace{-2.2cm}

\begin{center} \begin{LARGE} The Sine Transform of Isotropic Measures \\[0.5cm] \end{LARGE}

\begin{large}Gabriel Maresch and Franz E. Schuster \end{large}
\end{center}

\vspace{-0.8cm}

\begin{quote}
\footnotesize{ \vskip 1truecm\noindent {\bf Abstract.} Sharp
isoperimetric inequalities for the sine transform of even
isotropic measures are established. The corresponding reverse
inequalities are obtained in an asymptotically optimal form.
These new inequalities have direct applications to strong volume
estimates for convex bodies from data about their sections or
projections.}
\end{quote}

\vspace{0.4cm}

\centerline{\large{\bf{ \setcounter{abschnitt}{1}
\arabic{abschnitt}. Introduction}}} \alpheqn

\vspace{0.4cm}

The (spherical) cosine transform plays a fundamental role in
modern \linebreak geometric analysis. It arises naturally in a
number of different areas such as functional analysis, geometric
tomography and stochastic geometry (see e.g.,
\textbf{\cite{bourgain:lindenstrauss, gardner95, ludwig02,
Ludwig:Minkowski, Ryabogin:Zvavitch, schnweil, Thom96, YY06}}). A
classical theorem of Lewis \textbf{\cite{lewis78}} shows that
each finite dimensional subspace of $L_1$ is isometric to a
Banach space whose norm is the cosine transform of some even
isotropic measure on the unit sphere. This important result of
Lewis allows effortless proofs of isoperimetric inequalities,
which characterize Euclidean subspaces of $L_1$, by applications
of the Urysohn and H\"older inequalities (see
\textbf{\cite{LYZ2004}} for details). The reverse inequalities,
having $l_1^n$ subspaces and their duals as extremals, would turn
out to be significantly more difficult to establish. The
breakthrough here was achieved by Ball and Barthe in the 1990's.

Sharp reverse isoperimetric inequalities for the unit and polar
unit balls of subspaces of $L_1$ were established by Ball
\textbf{\cite{ball91a, ball91b}} using his ingenious \linebreak
reformulation of the Brascamp--Lieb inequality. The uniqueness
problem for the extremals in Ball's inequalities was solved by
Barthe \textbf{\cite{barthe98i}} for discrete isotropic measures
by using his newly obtained equality conditions for the
Brascamp--Lieb inequality and its inverse form. Recently, Lutwak,
Yang, and Zhang \textbf{\cite{LYZ2004, LYZ2007}} settled the
uniqueness questions for the extremal cases in Ball's
inequalities for general isotropic measures exploiting a more
direct approach based on the ideas of Ball and Barthe (see also
Barthe \textbf{\cite{barthe04}}).

In this article we obtain sharp isoperimetric inequalities for the
(spherical) \linebreak {\it sine transform} of isotropic measures.
We establish the corresponding reverse inequalities in an
asymptotically optimal form using the multi\-dimensional
Brascamp--Lieb inequality and its inverse obtained by Barthe
\textbf{\cite{barthe98i}}. While not as well known as the cosine
transform, its natural dual -- the sine transform -- appears in
different guises in geometric tomography. Therefore, applications
of our new inequalities lead to asymptotically sharp volume
estimates for convex bodies from certain data about their sections
or projections.

\pagebreak

The setting for this article is Euclidean $n$-space
$\mathbb{R}^n$ with $n \geq 3$. We use $\|\cdot\|$ to denote the
standard Euclidean norm on $\mathbb{R}^n$ and we write $x \cdot
y$ for the standard inner product of $x, y \in \mathbb{R}^n$. A
non-negative finite Borel measure $\mu$ on the unit sphere
$S^{n-1}$ is said to be {\it isotropic} if it has the same moment
of inertia about all lines through the origin or, equivalently,
if for all $x \in \mathbb{R}^n$,
\[\|x\|^2=\int_{S^{n-1}}|x \cdot u|^2\,d\mu(u).  \]

Isotropic measures have been the focus of recent studies, in
particular, in relation with a variety of extremal problems for
convex bodies (see, e.g., \textbf{\cite{giannomil, giamilrud,
giannopap, gruber08, LYZ2005}} and the references therein). Two
basic examples of isotropic measures on $S^{n-1}$ are (suitably
normalized) spherical Lebesgue measure and the cross measures,
i.e., measures concentrated uniformly on $\{\pm b_1, \ldots, \pm
b_n\}$, where $b_1,\ldots,b_n$ denote orthonormal basis vectors of
$\mathbb{R}^n$.

The {\it cosine transform} $\,\mathcal{C}\mu$ of a finite Borel
measure $\mu$ on $S^{n-1}$ is the \linebreak continuous function
defined by
\begin{equation} \label{cosin}
(\mathcal{C}\mu)(x)=\int_{S^{n-1}} |x \cdot u|\,d\mu(u), \qquad x
\in \mathbb{R}^n.
\end{equation}

If $\mu$ is not concentrated on a great subsphere and even (i.e.,
it assumes the same value on antipodal sets), its cosine
transform uniquely determines a norm on $\mathbb{R}^n$ whose unit
ball we denote by $C_{\mu}^*$.

In a highly influential paper, Bolker \textbf{\cite{bolker69}}
has shown that a convex body is the unit ball of an
$n$-dimensional subspace of $L_1$ if and only if the associated
norm admits a representation of the form (\ref{cosin}) for some
even measure $\mu$ not concentrated on a great subsphere.
Consequently, isoperimetric inequalities for the convex body
$C_{\mu}^*$ or its polar $C_{\mu}$ having ellipsoids as extremals
provide characterizations of Euclidean subspaces of $L_1$ (see
\textbf{\cite{LYZ2004}}).

Optimal reverse isoperimetric inequalities for the unit balls of
subspaces of $L_p$ -- having $l_p^n$ subspaces as extremals --
were established by Ball \textbf{\cite{ball91b}} using his
normalized Brascamp--Lieb inequality. The corresponding
inequalities for the polar unit balls of $L_1$ were also obtained
by Ball \textbf{\cite{ball91a}} and he predicted that for $p >
1$, these inequalities would follow from an inverse form of the
\linebreak Brascamp--Lieb inequality. Barthe
\textbf{\cite{barthe98i}} obtained this reverse Brascamp--Lieb
inequality in 1998 and used it to establish the reverse volume
inequalities for the polar unit balls of subspaces of $L_p$.
These landmark  results of Ball and Barthe have had a tremendous
impact on geometric analysis, see, e.g., \linebreak
\textbf{\cite{ball89, ball01, barthe98, bartcord04,
bartcordmaur06, tao08, carcor09, gardner02, LYZ2005, LYZ2009}}.

\pagebreak

The uniqueness questions for the extremal cases in the
inequalities of Ball and Barthe were completely settled only
recently by Lutwak, Yang, and Zhang \textbf{\cite{LYZ2004}} and
later, independently, by Barthe \textbf{\cite{barthe04}}. The
volume inequalities for subspaces of $L_1$ state the following:
Among even isotropic measures, $V(C_{\mu}^*)$ is maximized
precisely by cross measures and minimized precisely by normalized
Lebesgue measure, while $V(C_{\mu})$ is maximized precisely by
normalized Lebesgue measure and minimized precisely by cross
measures.

\vspace{0.3cm}

\noindent {\bf Definition} The \emph{sine transform}
$\,\mathcal{S}\mu$ of a finite Borel measure $\mu$ on $S^{n-1}$
is the continuous function defined by
\[(\mathcal{S}\mu)(x)=\int_{S^{n-1}} \|x|u^{\bot}\|\,d\mu(u), \qquad x \in \mathbb{R}^n.  \]

\vspace{0.2cm}

Here, $\|x|u^{\bot}\|$ is the length of the orthogonal projection
of $x$ onto the hyperplane orthogonal to $u$. If $\mu$ is even and
not concentrated on two anti\-podal points, its sine transform
uniquely determines (see Section 2 for details) a norm on
$\mathbb{R}^n$ whose unit ball we denote by $S_{\mu}^*$ and its
polar by $S_{\mu}$.

Let $\kappa_n$ denote the volume of the Euclidean unit ball in
$\mathbb{R}^n$ and define
\[\alpha_n :=  \frac{n(n-1)^{2n}}{\Gamma(n)^{1/(n-1)}}  \qquad \mbox{and}
\qquad \gamma_n := \frac{(n-1)\kappa_{n-1}^2}{\kappa_{n-2}\kappa_n}.  \]

The main results of this article are the following two theorems.

\begin{satz} \label{thm1} If $\mu$ is an even isotropic measure on $S^{n-1}$,
then
\begin{equation} \label{main1}
\frac{\kappa_n}{\gamma_n^n} \leq V(S_{\mu}^{\,*}) \leq
\frac{\kappa_n\gamma_n^n}{\alpha_n},
\end{equation}
with equality on the left if and only if $\mu$ is normalized
Lebesgue measure.
\end{satz}

While we believe that the right inequality in (\ref{main1}) is not
sharp for any value of $n$ (compare the discussion in Section 4),
we will show that it is {\it asymptotically optimal}. More
precisely, we will see in Section 4 that, up to a constant factor
which tends to one as $n$ goes to infinity, $V(S_{\mu}^*)$ is
maximized by cross measures.

\begin{satz} \label{thm2} If $\mu$ is an even isotropic measure on $S^{n-1}$,
then
\begin{equation} \label{main2}
\frac{\kappa_n\alpha_n}{\gamma_n^n} \leq V(S_{\mu}) \leq
\kappa_n\gamma_n^n,
\end{equation}
with equality on the right if and only if $\mu$ is normalized
Lebesgue measure.
\end{satz}

\pagebreak

We believe that the left inequality in (\ref{main2}) is also not
sharp. We will show, however, that it is asymptotically optimal.
More precisely, up to a factor \linebreak tending to one as $n$
goes to infinity, $V(S_{\mu})$ is minimized by cross measures.

The proofs of Theorems 1 and 2 are based on applications of the
Urysohn and H\"older inequalities, and the multidimensional
Brascamp--Lieb inequality and its inverse respectively. In our
approach we also make use of an instance of the Kantorovich
duality for the Brascamp--Lieb inequality and its inverse (see
Section 3 for details) which was exploited in their proof by
Barthe \textbf{\cite{barthe98i}}. This has the advantage that it
will provide additional geometric insight to the dual nature of
inequalities (\ref{main1}) and (\ref{main2}) (see Theorem 4.1).

The sine transform arises in geometric tomography in different
contexts (see Section 5 for a detailed account). In Section 5 we
will show that our main results -- Theorems 1 and 2 -- lead to
fairly strong volume estimates for convex bodies from certain
tomographic data which are dual to results of Giannopoulos and
Papadimitrakis \textbf{\cite{giannopap}} for the cosine transform.

\vspace{1cm}

\centerline{\large{\bf{ \setcounter{abschnitt}{2}
\arabic{abschnitt}. Background material}}}

\vspace{0.5cm} \reseteqn \alpheqn

For quick later reference, we collect in this section background
material regarding convex bodies. We also state some well known
facts about spherical harmonics which are needed to establish
basic injectivity properties of the sine transform. For a general
reference the reader may wish to consult the book by Schneider
\textbf{\cite{schneider93}}.

A convex body is a non-empty compact convex subset of
$\mathbb{R}^n$. We denote by $\mathcal{K}^n$ the space of convex
bodies in $\mathbb{R}^n$ endowed with the Hausdorff metric. A
convex body $K \in \mathcal{K}^n$ is uniquely determined by its
{\it support function} $h(K,\cdot)$, where $h(K,x)=\max \{x \cdot
y: y \in K\}$, $x \in \mathbb{R}^n$. Note that $h(K,\cdot)$ is
(positively) homogeneous of degree one and convex. Conversely,
each function with these properties is the support function of a
unique convex body.

The polar body $K^*$ of a convex body $K$ containing the origin
in its interior is defined by
\[K^*=\{x \in \mathbb{R}^n: x \cdot y \leq 1 \mbox{ for all } y \in K\}.  \]
Let $\rho(K,x)=\max\{\lambda \geq 0: \,\lambda\, x \in K\}$, $x
\in \mathbb{R}^n \backslash \{0\}$, denote the radial function of
$K$. It follows from the definitions of support functions and
radial functions, and the definition of the polar body of $K$,
that
\begin{equation} \label{suprad}
\rho(K^*,\cdot)=h(K,\cdot)^{-1} \qquad \mbox{and} \qquad
h(K^*,\cdot)=\rho(K,\cdot)^{-1}.
\end{equation}

\pagebreak

Using (\ref{suprad}) and the polar coordinate formula for volume,
it is easy to see that the volume of a convex body $K \in
\mathcal{K}^n$ containing the origin in its interior is given by
\begin{equation} \label{volexp}
V(K)=\frac{1}{n!}\int_{\mathbb{R}^n} \exp ( - h(K^*,x))\,dx,
\end{equation}
where integration is with respect to Lebesgue measure on
$\mathbb{R}^n$.

The classical Urysohn inequality (see, e.g.\ \textbf{\cite[\textnormal{p.\ 318}]{schneider93}}) provides an upper bound for the
volume of a convex body in terms of the average value of its support function: If $K \in \mathcal{K}^n$ has non-empty interior, then
\begin{equation} \label{urysohn}
\left (  \frac{V(K)}{\kappa_n} \right )^{1/n} \leq \frac{1}{n\kappa_n}\int_{S^{n-1}} h(K,u)\,du,
\end{equation}
with equality if and only if $K$ is a ball. Here the integral is
with respect to spherical Lebesgue measure.

For $K \in \mathcal{K}^n$ let $S(K)$ denote its surface area. The
well known classical isoperimetric inequality states that among
bodies of given volume, Euclidean balls have least surface area:
If $K \in \mathcal{K}^n$ has non-empty interior, then
\begin{equation} \label{iso}
n^n \kappa_nV(K)^{n-1} \leq S(K)^n,
\end{equation}
with equality if and only if $K$ is a ball.

Since convex bodies of a given volume may have arbitrarily large
surface area if they are very flat, the most natural way to
reverse the isoperimetric inequality is to consider affine
equivalence classes of convex bodies. This leads to the following
definition: The {\it minimal surface area} of a convex body $K \in
\mathcal{K}^n$ with non-empty interior is defined by
\[\partial(K) := \min\{ S(\phi K): \phi \in \mathrm{SL}(n)\}.  \]
We say $K$ is in {\it surface isotropic position} if $S(K) =
\partial(K)$. It was first proved by Petty \textbf{\cite{petty61}}
that every convex body with non-empty interior has a surface
isotropic position which is unique up to orthogonal
transformations.

The celebrated reverse isoperimetric inequality of Ball
\textbf{\cite{ball91b}} can now be stated as follows: If $K \in
\mathcal{K}^n$ has non-empty interior, then
\begin{equation} \label{reviso}
\partial(K)^n \leq
\frac{n^{3n/2}(n+1)^{(n+1)/2}}{n!}V(K)^{n-1}.
\end{equation}

\pagebreak

\noindent It was shown by Barthe \textbf{\cite{barthe98i}} that
equality holds in (\ref{reviso}) if and only if $K$ is a simplex.
It was also shown by Ball \textbf{\cite{ball91b}} that among
origin symmetric convex bodies of given volume the minimal
surface area is maximized (precisely) by the cube (the uniqueness
of extremals was settled by Barthe \textbf{\cite{barthe98i}}).

A convex body $K \in \mathcal{K}^n$ with non-empty interior is
also determined up to translations by its surface area measure
$S_{n-1}(K,\cdot)$. Recall that for a Borel set $\omega \subseteq
S^{n-1}$, $S_{n-1}(K,\omega)$ is the $(n-1)$-dimensional
Hausdorff measure of the set of all boundary points of $K$ at
which there exists a normal vector of $K$ belonging to $\omega$.
The following result of Petty \textbf{\cite{petty61}} (see also
\textbf{\cite{giannopap}}) will allow us to apply Theorems 1 and 2
in various geometric settings (see Section 5):

\begin{prop} \label{sip} A convex body $K \in \mathcal{K}^n$ with non-empty interior is in surface
isotropic position if and only if its surface area measure
$S_{n-1}(K,\cdot)$ is, up to normalization, isotropic.
\end{prop}

We conclude this section with a discussion of the injectivity
properties of the sine transform. To this end, we need some basic
facts about spherical harmonics (see e.g., Schneider
\textbf{\cite[\textnormal{Appendix}]{schneider93}}).

Let $\mathcal{H}^n_k$ denote the finite dimensional vector space
of spherical harmonics of dimension $n$ and order $k$ and let
$N(n,k)$ denote its dimension. We use \linebreak $L_2(S^{n-1})$ to
denote the Hilbert space of square integrable functions on
$S^{n-1}$ with its usual inner product $(\cdot\,,\cdot)$. The
spaces $\mathcal{H}^n_k$ are pairwise orthogonal with respect to
this inner product. In each space $\mathcal{H}_k^n$ we choose an
ortho\-normal basis $\{Y_{k1}, \ldots, Y_{kN(n,k)}\}$. Then
$\{Y_{k1}, \ldots, Y_{kN(n,k)}: k \in \mathbb{N}\}$ forms a
complete orthogonal system in $L_2(S^{n-1})$, i.e., for every $f
\in L_2(S^{n-1})$, the Fourier series
\[f \sim \sum \limits_{k=0}^{\infty}\mathrm{p}_k f  \]
converges in quadratic mean to $f$, where $\mathrm{p}_k f$ is the
orthogonal projection of $f$ onto $\mathcal{H}_k^n$. In
particular, for $f\in C(S^{n-1})$,
\begin{equation} \label{funcdet}
\mathrm{p}_kf=0 \quad \mbox{for all } k \in \mathbb{N} \qquad
\Longrightarrow \qquad f = 0.
\end{equation}
Thus, $f \in C(S^{n-1})$ is uniquely determined by its series
expansion.

For a finite Borel measure $\mu$ on $S^{n-1}$, we define
\[\mathrm{p}_k \mu=\sum\limits_{i=1}^{N(n,k)} \int_{S^{n-1}}Y_{ki}(u)\,d\mu(u)\, Y_{ki}.  \]
If $f \in C(S^{n-1})$, then
\[\left ( f,\mathrm{p}_k \mu \right ) = \int_{S^{n-1}} (\mathrm{p}_k f)(u)\,d\mu(u).   \]
Thus, by  (\ref{funcdet}), $\mu$ is uniquely determined by its
(formal) series expansion:
\begin{equation} \label{mudet}
\mathrm{p}_k \mu=0 \quad \mbox{for all } k \in \mathbb{N} \qquad
\Longrightarrow \qquad \mu = 0.
\end{equation}

A useful tool to establish injectivity results for integral
transforms is the Funk--Hecke theorem: Let $g$ be a continuous
function on $[-1,1]$. If $\mathrm{T}_{g}$ is the transformation
on the set of finite Borel measures on $S^{n-1}$ defined by
\begin{equation} \label{funkheck}
(\mathrm{T}_{g}\mu)(u)=\int_{S^{n-1}}g(u\cdot v)\,d\mu(v), \qquad
u \in S^{n-1},
\end{equation}
then there are real numbers $a_k[\mathrm{T}_{g}]$, the {\em
multipliers} of $\mathrm{T}_{g}$, such that
\[\mathrm{T}_{g}Y_k = a_k[\mathrm{T}_{g}]\,Y_k \]
for every $Y_k \in \mathcal{H}_k^n$. In particular, by Fubini's
theorem,
\begin{equation} \label{mult}
\mathrm{p}_k\! \left ( \mathrm{T}_{g}\mu \right
)=a_k[\mathrm{T}_{g}]\mathrm{p}_k \mu.
\end{equation}

Using (\ref{mudet}) and (\ref{mult}), it follows that a
transformation $\mathrm{T}_g$ defined on the space of finite Borel
measures on $S^{n-1}$ and satisfying (\ref{mult}) is injective if
and only if all multipliers $a_k[\mathrm{T}_g]$ are non-zero.

Clearly, the sine transform $\mathcal{S}$ of finite Borel measures
on $S^{n-1}$ is of the form (\ref{funkheck}), where
\[g(t)=\sqrt{1-t^2}, \qquad t \in [-1,1].  \]
Thus, by the Funk--Hecke theorem, the sine transform satisfies
(\ref{mult}). The multipliers $a_k[\mathcal{S}]$ have been
calculated in \textbf{\cite{goodeyweil92}}: For every $k \in
\mathbb{N}$, we have
\begin{equation} \label{sinmult}
a_{2k}[\mathcal{S}] \neq 0 \qquad \mbox{and} \qquad
a_{2k+1}[\mathcal{S}] = 0.
\end{equation}

Since $f \in C(S^{n-1})$ (or a measure $\mu$ on $S^{n-1}$) is even
if and only if $\mathrm{p}_k f=0$ (or $\mathrm{p}_k \mu$ = 0,
respectively) for every odd $k \in \mathbb{N}$, (\ref{sinmult})
yields the following injectivity result (for a stability version
see \textbf{\cite{hugschneider}}):

\begin{prop} \label{injsine} The sine transform is injective on even measures on
$S^{n-1}$.
\end{prop}

\pagebreak

\centerline{\large{\bf{ \setcounter{abschnitt}{3}
\arabic{abschnitt}. The Brascamp--Lieb inequality and its
inverse}}}

\vspace{0.5cm} \reseteqn \alpheqn \setcounter{koro}{0}

In the following we recall the rank $n - 1$ case of the
multidimensional Brascamp--Lieb inequality and its reverse form
which are crucial in the proof of our main results. We also state
a duality formula for these inequalities \linebreak established by
Barthe \textbf{\cite{barthe98i}} which is a special case of the
Kantorovich duality principle from optimal mass transportation
(see e.g.\ \textbf{\cite[\textnormal{Chapters 1 \&
6}]{villani03}}).

The Brascamp--Lieb inequality \textbf{\cite{braslieb, lieb}} was
established to prove the sharp form of Young's convolution
inequality. It's multidimensional version unifies and generalizes
several fundamental inequalities from geometric analysis such as
the H\"older inequality and the Loomis--Whitney inequality.

Around 1990 Ball \textbf{\cite{ball89}} discovered an important
reformulation of the Brascamp--Lieb inequality (later generalized
by Barthe \textbf{\cite{barthe98i}}) which exploited an additional
geometric hypothesis of the given data and was tailor-made for
applications in convex geometry. This geometric Brascamp--Lieb
inequality also allowed a simple computation of the optimal
constant. We will only need (and thus only state) this powerful
inequality in the rank $n - 1$ case.

In the following we write $\pi_u$, $u \in S^{n-1}$, for the
orthogonal projection onto the hyperplane $u^{\bot}$.

\vspace{0.4cm}

\noindent {\bf The Brascamp--Lieb Inequality. (\!\!\cite{lieb})}
\emph{Let $u_1,\ldots, u_m \in S^{n-1}$, $m \geq n$, and $c_1,
\ldots, c_m
> 0$ such that
\[\sum_{i=1}^m c_i \pi_{u_i}=\mathrm{Id}.   \]
If $f_i: u_i^{\bot} \rightarrow [0,\infty)$, $1 \leq i \leq m$,
are integrable functions, then
\begin{equation} \label{bralieb}
\int_{\mathbb{R}^n} \prod_{i=1}^m f_i(x|u_i^{\bot})^{c_i}dx \leq
\prod_{i=1}^m \left ( \int_{u_i^{\bot}} f_i  \right )^{c_i}.
\end{equation}
There is equality if the $f_i$, $1 \leq i \leq m$, are identical
Gaussian densities.}

\vspace{0.4cm}

The problem of characterizing all extremizers for the
multidimensional Brascamp--Lieb inequality was settled only
recently by Valdimarsson \textbf{\cite{valdi08}} after previous
contributions by a number of mathematicians (see
\textbf{\cite{barthe98i, tao08, carlielos04}}). \linebreak In
order to discuss the quality of our upper bound in Theorem 1, and
our lower bound in Theorem 2 respectively, we state the following
special case of this characterization for the rank $n - 1$ case
(cf.\ \textbf{\cite[\textnormal{Theorem 12}]{valdi08}}):

\pagebreak

\begin{prop} \label{orthoprop} Let $u_i \in S^{n-1}$,  $c_i > 0$, $1 \leq i \leq
m$, be as above. Suppose that $f_i: u_i^{\bot} \rightarrow
[0,\infty)$, $1 \leq i \leq m$, are (non-identically-zero)
integrable functions such that none of them is a Gaussian. If
equality holds in (\ref{bralieb}), then there exist an
orthonormal basis $\{b_1,\ldots,b_n\}$ of $\mathbb{R}^n$,
integrable functions $\varphi_i$ of one variable and constants
$a_i \in \mathbb{R}$, $1 \leq i \leq m$, such that
\[\{u_1,\ldots,u_m\} \subseteq \{\pm b_1, \ldots,\pm b_n\}\]
 and
\[f_i(x_1,\ldots,x_{i-1},x_{i+1},\ldots,x_n)=a_i\varphi_1(x_1)\cdots\varphi_{i-1}(x_{i-1})\varphi_{i+1}(x_{i+1})\cdots \varphi_n(x_n).   \]
\end{prop}

\vspace{0.1cm}

The strength of the Brascamp--Lieb inequality for volume
estimates of sections of the unit ball of $l_p^n$ was exploited
by Ball (see \textbf{\cite{ball89, ball91b}}). He also predicted
that a reverse form of the Brascamp--Lieb inequality would lead to
dual estimates for projections of the unit ball of $l_p^n$. The
breakthrough here was achieved by Barthe \textbf{\cite{barthe97,
barthe98i}} who established the reverse Brascamp--Lieb
inequality. In the following we state this inequality in the rank
$n - 1$ case which is needed in the proof of Theorem 2.

\vspace{0.5cm}

\noindent {\bf The Reverse Brascamp--Lieb Inequality.
(\!\!\cite{barthe98i})} \emph{Let $u_1,\ldots, u_m \in S^{n-1}$,
$m \geq n$, and $c_1, \ldots, c_m
> 0$ such that
\[\sum_{i=1}^m c_i \pi_{u_i}=\mathrm{Id}.   \]
If $f_i: u_i^{\bot} \rightarrow [0,\infty)$, $1 \leq i \leq m$,
are integrable functions, then
\begin{equation} \label{revbralieb}
\int_{\mathbb{R}^n}\! \sup \left \{ \prod_{i=1}^m f_i(y_i)^{c_i}:
x=\sum_{i=1}^m c_iy_i,\, y_i \in u_i^{\bot} \right \} dx \geq
\prod_{i=1}^m \left ( \int_{u_i^{\bot}}\! f_i \right )^{c_i}\!.
\end{equation}
There is equality if the $f_i$, $1 \leq i \leq m$, are identical
Gaussian densities.}

\vspace{0.5cm}

The proof of the reverse Brascamp--Lieb inequality by Barthe
relies on the \linebreak existence and uniqueness of a certain
measure preserving map, the socalled \linebreak Brenier map,
between two sufficiently regular probability measures (see e.g.\
\textbf{\cite{brenier91, McCann95}}). Barthe's proof also
exploited a classical principle dating back to Kantorovich which
states that the problem of optimal mass transportation admits two
dual formulations. In particular, this duality principle made it
possible to derive both the Brascamp--Lieb inequality and its
inverse from a single inequality which is stated in the following
theorem.

\pagebreak

\begin{theorem} \label{kanttheo} Let
$u_1,\ldots, u_m \in S^{n-1}$, $m \geq n$, and $c_1, \ldots, c_m
> 0$ such that
\[\sum_{i=1}^m c_i \pi_{u_i}=\mathrm{Id}.   \]
If $f_i, g_i: u_i^{\bot} \rightarrow [0,\infty)$, $1 \leq i \leq
m$, are integrable functions such that
\[ \int_{u_i^{\bot}} f_i = \int_{u_i^{\bot}} g_i = 1, \]
then
\begin{equation} \label{kantdual}
\int_{\mathbb{R}^n} \prod_{i=1}^m f_i(x|u_i^{\bot})^{c_i}dx \leq\!
\int_{\mathbb{R}^n}\! \sup \left \{ \prod_{i=1}^m g_i(y_i)^{c_i}:
x=\sum_{i=1}^m c_iy_i, y_i \in u_i^{\bot} \right \} dx.
\end{equation}
\end{theorem}

\vspace{0.2cm}

Note that equality in (\ref{kantdual}) can only hold if the $f_i$
are extremizers for the Brascamp--Lieb inequality and the $g_i$
are extremizers for the reverse Brascamp--Lieb inequality.

Inequality (\ref{kantdual}) will provide a convenient way to
obtain the upper bound in Theorem 1 and the lower bound in Theorem
2 from a single inequality.

\vspace{1cm}

\centerline{\large{\bf{ \setcounter{abschnitt}{4}
\arabic{abschnitt}. Proof of the main results}}}

\vspace{0.5cm} \reseteqn \alpheqn \setcounter{koro}{0}

After these preparations, we are now in a position to prove our
main theorems. In fact we will establish stronger results since we
consider in this section arbitrary (and not necessarily even)
isotropic measures.

The following two results, which directly imply Theorems 1 and 2,
make use of Theorem \ref{kanttheo} in our context:

\begin{theorem}\label{satzupper}
If $\mu$ is an isotropic measure on $S^{n-1}$, then
\[V(S_{\mu}^*) \leq V(S_{\mu})/\alpha_n.  \]
\end{theorem}

\noindent {\it Proof}\,: First assume that $\mu$ is discrete, say
$\mathrm{supp}\, \mu = \{u_1, \ldots, u_m\}$ and $\mu(\{u_i\})=:
\bar{c}_i > 0$. Since $\mu$ is isotropic, it follows that
$\mu(S^{n-1})=\sum_{i=1}^m\bar{c}_i=n$. Therefore, using $\pi_u =
\mathrm{Id}-u\otimes u$, we have
\begin{equation} \label{zerl17}
\frac{1}{n-1}\sum_{i=1}^m \bar{c}_i\, \pi_{u_i} = \sum_{i=1}^m
\bar{c}_i\, u_i \otimes u_i =\mathrm{Id}.
\end{equation}
From (\ref{volexp}) and the definition of the sine transform, it
follows that
\begin{equation} \label{vol171}
V(S_{\mu}^*)=\frac{1}{n!}\int_{\mathbb{R}^n}\prod_{i=1}^m\exp(-(n-1)\|x|u_i^{\bot}\|)^{c_i}dx,
\end{equation}
where $c_i:=\bar{c}_i/(n-1)$, $i = 1, \ldots, m$. Let $B$ denote
the Euclidean unit ball in $\mathbb{R}^n$. Since
$\|x|u^{\bot}\|=h(B|u^{\bot},x)$, we have
\[S_{\mu} = \left \{x \in \mathbb{R}^n: x=\sum_{i=1}^m \bar{c}_i y_i,\,y_i \in B|u_i^{\bot}\right \}.  \]
Consequently, we obtain
\begin{equation} \label{vol172}
V(S_{\mu}) = \int_{\mathbb{R}^n} \sup\left \{ \prod_{i=1}^m
\mathbf{1}_{[0,n-1]}(\|y_i\|)^{c_i}: x=\sum_{i=1}^m c_i y_i,\,y_i
\in u_i^{\bot} \right \}dx.
\end{equation}
Define functions $f_i,g_i: u_i^{\bot} \rightarrow [0,\infty), 1
\leq i \leq m$, by
\begin{equation} \label{deffi}
f_i(y)=\frac{(n-1)^{n-1}}{\Gamma(n)\kappa_{n-1}}\exp(-(n-1)\|y\|)
\end{equation}
and
\begin{equation} \label{defgi}
g_i(y)=\frac{1}{(n-1)^{n-1}\kappa_{n-1}}\mathbf{1}_{[0,n-1]}(\|y\|).
\end{equation}
Note that the normalizations are chosen such that
\[\int_{u_i^{\bot}} f_i = \int_{u_i^{\bot}} g_i = 1.  \]
Since $\sum_{i=1}^mc_i=n/(n-1)$, we obtain, by (\ref{zerl17}) --
(\ref{vol172}) and Theorem \ref{kanttheo},
\begin{align} \label{inequ217}
V(S_{\mu}^*) & = \frac{(\Gamma(n)\kappa_{n-1})^{n/(n-1)}}{n!(n-1)^n}\!\int_{\mathbb{R}^n}\prod_{i=1}^mf_i(x|u_i^{\bot})^{c_i}dx \nonumber \\
& \leq \frac{(\Gamma(n)\kappa_{n-1})^{n/(n-1)}}{n!(n-1)^n}\!
\int_{\mathbb{R}^n}\! \sup \left \{ \prod_{i=1}^m g_i(y_i)^{c_i}:
x=\sum_{i=1}^m c_iy_i, y_i \in u_i^{\bot} \right \} dx \\
& =  \frac{\Gamma(n)^{1/(n-1)}}{n(n-1)^{2n}}V(S_{\mu})=
V(S_{\mu})/\alpha_n. \nonumber
\end{align}

Now let $\mu$ be an arbitrary isotropic measure on $S^{n-1}$. As
in \textbf{\cite[\textnormal{pp.\ 55--56}]{barthe04}} construct a
sequence $\mu_k$, $k \in \mathbb{N}$, of discrete isotropic
measures such that $\mu_k$ converges weakly to $\mu$ as $k
\rightarrow \infty$. It follows that $\lim_{k \rightarrow \infty}
h(S_{\mu_k},v) = h(S_{\mu},v)$ for every $v \in S^{n-1}$. Since
the pointwise convergence of support functions implies the
convergence of the respective convex bodies in the Hausdorff
metric (see, e.g., \textbf{\cite[\textnormal{Chapter
1}]{schneider93}}), the continuity of volume and polarity on
convex bodies containing the origin in their interiors finishes
the proof. \hfill $\blacksquare$

\vspace{0.4cm}

Our next result completes the proof of Theorems 1 and 2:

\begin{theorem} \label{satzlower} If $\mu$ is an isotropic measure on
$S^{n-1}$, then
\[\frac{\kappa_n}{\gamma_n^n} \leq V(S_{\mu}^{\,*}) \qquad \mbox{and} \qquad V(S_\mu) \leq \kappa_n \gamma_n^n.  \]
If $\mu$ is even, then there is equality in either inequality if
and only if $\mu$ is normalized Lebesgue measure.
\end{theorem}

\noindent {\it Proof}\,: By the polar coordinate formula for
volume, (\ref{suprad}), and the H\"older inequality, we have
\[\left( \frac{V(S_{\mu}^*)}{\kappa_n}\right)^{-1/n}=\left(
\frac{1}{n\kappa_n} \int_{S^{n-1}} h(S_{\mu},u)^{-n}\, du
\right)^{-1/n} \le \frac{1}{n\kappa_n}\int_{S^{n-1}}
h(S_{\mu},u)\, du \] with equality if and only if
$h(S_{\mu},\cdot)$ is constant, i.e.\ $S_{\mu}$ is a ball. From
the definition of the sine transform and Fubini's theorem, we
obtain
\begin{align*}
\frac{1}{n\kappa_n}\int_{S^{n-1}} h(S_{\mu},u)\, du & =
\frac{1}{n\kappa_n} \int_{S^{n-1}} \int_{S^{n-1}} \sqrt{1-(u
\cdot v)^2}\, du\, d\mu(v) \\
& = (n-1)\frac{\kappa_{n-1}}{\kappa_n}\int_{-1}^1 (1-t^2)^{n/2-1}
dt = \gamma_n.
\end{align*}
Consequently,
\[\left( \frac{V(S_{\mu}^*)}{\kappa_n}\right)^{-1/n} \leq \gamma_n  \]
with equality if and only if $S_{\mu}$ is a ball. Proposition
\ref{injsine} now yields the equality conditions for even
isotropic measures.

In order to establish the second inequality, we apply the
classical Urysohn inequality (\ref{urysohn}) to obtain
\[\left( \frac{\vol(S_{\mu})}{\kappa_n}\right)^{1/n} \le \frac{1}{n\kappa_n} \int_{S^{n-1}} h(S_{\mu},u)\, du = \gamma_n \]
with equality if and only if $S_{\mu}$ is a ball. Again, the
equality conditions for even isotropic measures follow from
Proposition \ref{injsine}. \hfill $\blacksquare$

\vspace{0.4cm}

We do not believe that our upper bound in Theorem \ref{thm1} and
our lower bound in Theorem \ref{thm2} are sharp: For equality to
hold in these inequalities, we must have equality in Theorem
\ref{satzupper}. If the isotropic measure $\mu$ is discrete, this
is equivalent to equality in (\ref{inequ217}). But from the remark
after Theorem \ref{kanttheo}, \linebreak Proposition
\ref{orthoprop} and the specific form of the functions $f_i$,
$g_i$ defined in (\ref{deffi}) and (\ref{defgi}), it follows that
equality can not hold in (\ref{inequ217}). Hence, for discrete
measures we can not have equality in Theorem \ref{satzupper}. In
order to deduce the same fact for arbitrary isotropic measures,
we need a continuous analogue of Proposition \ref{orthoprop}.
Such a result was proved by Barthe
\textbf{\cite[\textnormal{Theorem 2}]{barthe04}} for the rank 1
case of the Brascamp--Lieb inequality using the equality
conditions for a determinant inequality by Ball which were
obtained by Lutwak, Yang, and Zhang \textbf{\cite[\textnormal{p.\
168}]{barthe04}}. Unfortunately, neither Barthe's continuous
analogue of Proposition \ref{orthoprop} nor the equality
conditions of Ball's determinant inequality are known in (the more
complex) rank $n - 1$ case.

The following result shows, however, that our upper bound in
Theorem \ref{thm1} \linebreak and our lower bound in Theorem
\ref{thm2}, respectively, are asymptotically optimal in a strong
sense:

\begin{theorem} \label{asopt} If $\nu_n$, $n \geq 3$, are cross measures on
$S^{n-1}$, then
\[ \lim_{n \rightarrow \infty} \frac{\alpha_n}{\kappa_n\gamma_n^n} V(S_{\nu_n}^{\,*})
=\lim_{n \rightarrow \infty} \frac{\gamma_n^n}{\kappa_n\alpha_n}
V(S_{\nu_n})=1.  \]
\end{theorem}

\noindent {\it Proof}\,: Let $\mathrm{supp}\,\nu_n = \{\pm e_1,
\ldots,\pm e_n\}$, where $\{e_1, \ldots, e_n\}$ is an orthonormal
basis of $\mathbb{R}^n$. By the definition of the sine transform,
the support function of $S_{\nu_n}$ is given by
\[h(S_{\nu_n},v) = \sum_{i=1}^n \|\pi_{e_i} v\|, \qquad v \in S^{n-1}.   \]

A simple computation shows that
\[\max_{v \in S^{n-1}}\sum_{i=1}^n \|\pi_{e_i} v\| = n\sqrt{1-\frac 1n}.   \]
(The maximum is attained precisely at the points $(\pm
\frac{1}{\sqrt{n}}, \ldots, \pm \frac{1}{\sqrt{n}})$.) Hence, we
have the inclusion
\begin{equation} \label{incl17}
S_{\nu_n} \subseteq n\sqrt{1-\frac 1n}\, B.
\end{equation}

Theorem \ref{thm1} and Theorem \ref{thm2} together with
(\ref{incl17}) now immediately yield the following volume bounds
for $S_{\nu_n}^*$ and $S_{\nu_n}$, respectively:
\begin{equation} \label{stirr1}
\frac{\kappa_n}{n^n} \left(1-\frac 1n\right)^{-n/2} \leq
V(S_{\nu_n}^*) \leq \frac{\kappa_n\gamma_n^n}{\alpha_n} \quad
\qquad \qquad \phantom{a}
\end{equation}
and
\begin{equation} \label{stirr2}
\phantom {a} \quad \qquad \qquad \frac{\kappa_n
\alpha_n}{\gamma_n^n} \leq V(S_{\nu_n}) \le \kappa_n n^n
\left(1-\frac 1n\right)^{n/2}.
\end{equation}

Using Stirling's formula and the definition of the constants
$\alpha_n$ and $\gamma_n$, it is easy to show that
\[\lim_{n\rightarrow \infty} \frac{\alpha_n}{n^n\gamma_n^n}\left(1-\frac
1n\right)^{-n/2} = 1.  \] Consequently, we also have
\[\lim_{n\rightarrow \infty} \frac{n^n\gamma_n^n}{\alpha_n}\left(1-\frac
1n\right)^{n/2} = 1   \] which completes the proof in view of
(\ref{stirr1}) and (\ref{stirr2}). \hfill $\blacksquare$

\vspace{0.4cm}

In view of Theorem \ref{asopt}, we formulate the following

\vspace{0.4cm}

\noindent {\bf Conjecture.} \emph{Among even isotropic measures,
$V(S_{\mu}^*)$ is maximized precisely by cross measures, while
$V(S_{\mu})$ is minimized precisely by cross measures.}

\vspace{1cm}

\centerline{\large{\bf{ \setcounter{abschnitt}{5}
\arabic{abschnitt}. The sine transform in geometric tomography}}}
\alpheqn

\vspace{0.5cm} \reseteqn \alpheqn \setcounter{koro}{0}

In this last section we briefly recall several tomographic
operators on convex bodies induced by the sine transform. As
applications of our main results, we then present asymptotically
optimal volume inequalities for these operators. Our results are
dual to volume estimates due to Giannopoulos and Papadimitrakis
\textbf{\cite{giannopap}} for projection bodies which we also
recall.

The {\it projection body} $\Pi K$ of $K \in \mathcal{K}^n$ is the
convex body defined by
\[h(\Pi K,v)=\mathrm{vol}_{n-1}(K|v^{\bot})=\frac{1}{2}\int_{S^{n-1}}|u\cdot v|\,dS_{n-1}(K,v), \qquad v \in S^{n-1}.  \]
Here, the second equation is the well known Cauchy projection
formula.

\pagebreak

Projection bodies were introduced by Minkowski at the turn of the
\linebreak previous century and have since become an important
tool in the study of projections of convex bodies (see e.g.\
\textbf{\cite{ball91a, bolker69, bourgain:lindenstrauss,
kiderlen08, Ryabogin:Zvavitch, Schu06a, Schu08, Schu10}}). It was
first proved by Petty \textbf{\cite{petty67}} that for all $\phi
\in \mathrm{GL}(n)$ and all $K \in \mathcal{K}^n$,
\begin{equation} \label{pigln}
\Pi(\phi K)=|\det \phi\,|\phi^{-\mathrm{T}}\,\Pi K.
\end{equation}

In particular, (\ref{pigln}) shows that the volume of projection
bodies and their polars is invariant under volume preserving
linear transformations. The \linebreak fundamental affine
isoperimetric inequalities for polar projection bodies are the
Petty \textbf{\cite{petty72}} and the Zhang
\textbf{\cite{zhang91}} projection inequalities (for an important
recent generalization of Petty's projection inequality, see
\textbf{\cite{LYZ2000a}}): If $K \in \mathcal{K}^n$ has non-empty
interior, then
\begin{equation} \label{petzhang}
\frac{(2n)!}{n^n(n!)^2} \leq V(\Pi^* K)V(K)^{n-1} \leq \left (
\frac{\kappa_n}{\kappa_{n-1}}  \right )^n.
\end{equation}
There is equality in the left inequality if and only if $K$ is a
simplex and equality in the right inequality if and only if $K$
is an ellipsoid. It is a major open problem to determine the
corresponding inequalities for the volume of the projection body
itself.

In \textbf{\cite{giannopap}} Giannopoulos and Papadimitrakis
first observed that the volume inequalities of Ball for unit and
polar unit balls of subspaces of $L_1$ admit an elegant
reformulation using projection bodies (cf.\ the proof of Theorem
\ref{singeo}):

\begin{theorem} \label{giapap} \textbf{\emph{(\!\cite{giannopap})}} If $K \in \mathcal{K}^n$ has non-empty interior, then
\[\kappa_n\left (\frac{n\kappa_n}{\kappa_{n-1}}\right )^n\leq V(\Pi^* K)\partial(K)^n \leq \frac{4^nn^n}{n!} \quad \qquad \phantom{a}   \]
and
\begin{equation} \label{proj17}
\phantom{a\,\,} \quad \qquad \frac{1}{n^n} \leq V(\Pi
K)/\partial(K)^{n} \leq \left ( \frac{\kappa_{n-1}}{n\kappa_n}
\right )^n \kappa_n.
\end{equation}
\end{theorem}

Note that all the inequalities of Theorem \ref{giapap} are sharp;
consider e.g.\ \linebreak ellipsoids and parallelotopes. For
centrally-symmetric convex bodies the equality conditions were
settled by Lutwak, Yang, and Zhang in \textbf{\cite{LYZ2004}}.

The inequalities (\ref{proj17}) together with the isoperimetric
inequality (\ref{iso}) and its exact reverse form (\ref{reviso})
immediately provide asymptotically optimal \linebreak reverse
forms of the Petty and Zhang projection inequalities
(\ref{petzhang}):

\begin{koro} \label{openppi} If $K \in \mathcal{K}^n$ has non-empty interior, then
\[ n^{-1/2} \leq \left [ V(\Pi K)/V(K)^{n-1} \right ]^{1/n} \leq e^{3/2}.   \]
\end{koro}

\vspace{0.2cm}

Up to a constant multiple, both inequalities in Corollary
\ref{openppi} are best possible; consider e.g.\ ellipsoids and
simplices (see also \textbf{\cite{LYZ2001}}).

\vspace{0.2cm}

The sine transform of surface area measures also arises naturally
in \linebreak geometric tomography in a number of different
guises.

\vspace{0.4cm}

\noindent {\bf Examples:}
\begin{enumerate}
\item[(a)] If $K \in \mathcal{K}^n$ and $v \in S^{n-1}$, then it
was shown by Schneider \textbf{\cite{schn70}} that
\[\int_{-\infty}^{\infty}V_{n-2}(K \cap (v^{\bot} + tv))\,dt = \frac{1}{2(n+1)} \int_{S^{n-1}} \|v|u^{\bot}\|\,dS_{n-1}(K,u),  \]
where $2V_{n-2}(L)$ denotes the $(n-2)$-dimensional surface area
of an $(n-1)$-dimensional convex body $L$. Thus, the sine
transform of the surface area measure of $K$ is, up to a factor,
the integrated surface area of parallel hyperplane sections of
$K$.

\item[(b)] For $i \in \{1, \ldots, n - 1\}$, the $i$-th
mean section operator $\mathrm{M}_i: \mathcal{K}^n \rightarrow
\mathcal{K}^n$, introduced by Goodey and Weil in
\textbf{\cite{goodeyweil92}}, is defined by
\[h(\mathrm{M}_iK,\cdot) = \int_{\mathrm{AGr}_{i,n}} h(K \cap E, \cdot)\,d\sigma_i(E).  \]
Here, $\mathrm{AGr}_{i,n}$ is the affine Grassmannian of
$i$-dimensional planes in $\mathbb{R}^n$ and $\sigma_i$ is its
(suitably normalized) motion invariant measure. It was shown in
\textbf{\cite{goodeyweil92}} that for origin-symmetric convex
bodies
\[h(\mathrm{M}_2K,\cdot) =  \frac{\kappa_2^2 \kappa_{n-2}}{n(n-1)\kappa_n}\int_{S^{n-1}} \|v|u^{\bot}\|\,dS_{n-1}(K,u).   \]
\item[(c)] For $i \in \{0, \ldots, n\}$, let $V_i(K)$ denote the $i$-th
intrinsic volume of $K \in \mathcal{K}^n$. The projection body
$\Pi_i K$ of order $i$ of $K$ is defined by
\[h(\Pi_i K,v)=V_i(K|v^{\bot}), \qquad v \in S^{n-1}.\]
A direct computation shows that
\[h(\Pi_1 \Pi K,v)=\frac{\kappa_{n-2}}{n-1}\int_{S^{n-1}}
\|v|u^{\bot}\|\,dS_{n-1}(K,u).\]
\end{enumerate}

\pagebreak

An important part of geometric tomography deals with the
estimation of the volume (and other quantities) of a convex or
star body from data about the projections or the sections of the
body (see e.g.\ \textbf{\cite{ball91a, Gardner99, haberl09,
Ryabogin:Zvavitch, YY06, wannerer10}} \linebreak and, in
particular, \textbf{\cite[\textnormal{Chapter 9}]{gardner95}} and
the references therein). The rest of this section is devoted to
establishing volume inequalities for the examples above which are
(in some sense) dual to Theorem \ref{giapap} and Corollary
\ref{openppi}. In order to allow for an immediate comparison with
the results for projection bodies, we will introduce yet another
operator $\Psi: \mathcal{K}^n \rightarrow \mathcal{K}^n$, defined
by
\[h(\Psi K,v) = \frac{\kappa_{n-2}}{(n-1)\kappa_{n-1}} \int_{S^{n-1}} \|v|u^{\bot}\|\,dS_{n-1}(K,u).\]
Here, the normalization is chosen such that $\Pi B = \Psi B$.

It is important to note that while $\Psi$ still commutes with
orthogonal transformations, it does {\it not} intertwine affine
transformations like the projection body map $\Pi$. (The very
special role of the projection body operator in affine convex
geometry has only been demonstrated recently by Ludwig
\textbf{\cite{ludwig02, Ludwig:Minkowski}}.) Consequently, the
quantities $V(\Psi K)$ and $V(\Psi^*K)$ are rigid motion
invariant but not invariant under volume preserving linear
transformations. In fact, for a convex body $K$ of given volume,
$V(\Psi K)$ may be arbitrarily large and $V(\Psi^* K)$
arbitrarily small, respectively. We will therefore fix a position
of the body, to be more precise, the surface isotropic position,
to bound the quantities $V(\Psi K)$ and $V(\Psi^* K)$.

The following result is a reformulation of the slightly more
general versions of Theorems \ref{thm1} and \ref{thm2} proved in
Section 4:

\begin{theorem} \label{singeo} If $K \in \mathcal{K}^n$ is in surface isotropic
position, then
\begin{equation*}
\phantom{a} \,\,\, \qquad \qquad \kappa_n\left
(\frac{n\kappa_n}{\kappa_{n-1}}\right )^n \leq V(\Psi^*
K)\partial(K)^n \leq \left ( \frac{n}{\kappa_{n}} \right )^{n-1}
\frac{\kappa_{n-1}^{3n} \Gamma(n)^{1/(n-1)}}{\kappa_{n-2}^{2n}},
\end{equation*}
with equality among centrally-symmetric convex bodies in the left
inequality if and only if $K$ is a ball, and
\begin{equation*}
\frac{\kappa_{n-2}^{2n}\kappa_n^2}{\kappa_{n-1}^{3n}
\Gamma(n)^{1/(n-1)}} \left ( \frac{\kappa_n}{n} \right )^{n-1}
\leq V(\Psi K)/
\partial(K)^n \leq \left ( \frac{\kappa_{n-1}}{n\kappa_n} \right
)^n \kappa_n, \qquad \qquad \quad \phantom{a}
\end{equation*}
with equality among centrally-symmetric convex bodies in the
right inequality if and only if $K$ is a ball.
\end{theorem}

\noindent {\it Proof}\,: Define the non-negative Borel measure
$\mu$ on $S^{n-1}$ by
\[\mu = \frac{n}{\partial(K)}S_{n-1}(K,\cdot).  \]
Since $K$ is in surface isotropic position, it follows from
Proposition \ref{sip} that $\mu$ is isotropic. Clearly, by the
definitions of $S_{\mu}$ and the map $\Psi$, we have
\[S_{\mu}^* = \frac{\partial(K)\kappa_{n-2}}{n(n-1)\kappa_{n-1}}\Psi^* K \qquad \mbox{and} \qquad S_{\mu} = \frac{n(n-1)\kappa_{n-1}}{\partial(K)\kappa_{n-2}}\Psi K .  \]
Applications of Theorem \ref{satzupper} and Theorem
\ref{satzlower}, now complete the proof. \hfill $\blacksquare$

\vspace{0.4cm}

A combination of the inequalities of Theorem \ref{singeo} with the
isoperimetric inequality (\ref{iso}) and its exact reverse form
(\ref{reviso}), now yields

\begin{koro} \label{sinass} If $K \in \mathcal{K}^n$ is in surface isotropic position, then
\[\phantom{a} \,\,\,\, (en)^{-1} \leq \left [ V(\Psi^* K)V(K)^{n-1} \right ]^{1/n} \leq e^{3/2}n^{-1/2} \]
and
\[ n^{-1/2} \leq \left [ V(\Psi K)/V(K)^{n-1} \right ]^{1/n} \leq e^{3/2}. \,\,\,\, \phantom{a}  \]
\end{koro}

\vspace{0.2cm}

Note that both bounds are, up to a constant multiple, best
possible \linebreak (consider e.g.\ Euclidean balls and cubes) and
that they are precisely of the same order as the corresponding
bounds for projection bodies given by (\ref{petzhang}) and
Corollary \ref{openppi}.

We finally remark that it is easy to show that among convex
bodies of given volume there exists an upper bound for the
quantity $V(\Psi^*K)$, and a lower bound for $V(\Psi K)$
respectively. It is the authors believe that both bounds are
attained precisely by Euclidean balls. For convex bodies in
surface isotropic position, Corollary \ref{sinass} confirms these
conjectures asymptotically.

\vspace{0.5cm}

\medskip\noindent{\bf Acknowledgements.}
The authors are grateful to Keith Ball and Franck Barthe for very
helpful discussions on early versions of this manuscript. The work
of the second author was supported by the Austrian Science Fund,
within the project ``Minkowski valuations and geometric
inequalities", Project Number: P\,22388-N13.

\vspace{0.5cm}

Vienna University of Technology \par Institute of Discrete
Mathematics and Geometry \par Wiedner Hauptstra\ss e 8--10/1046
\par A--1040 Vienna, Austria \par
gabriel.maresch@tuwien.ac.at \par franz.schuster@tuwien.ac.at

\end{document}